\documentclass{article}
\usepackage{color}
\usepackage{amsmath, amssymb, theorem, latexsym}
\usepackage{graphicx}
\usepackage{epsf}

\theoremstyle{plain} %% Pierre
\newtheorem{theorem}{Theorem}[section]
\newtheorem{lemma}[theorem]{Lemma}
\newtheorem{proposition}[theorem]{Proposition}
\newtheorem{conjecture}[theorem]{Conjecture}
\newtheorem{corollary}[theorem]{Corollary}

\newtheorem{property}[theorem]{Property}
\newtheorem{properties}[theorem]{Properties}

\newtheorem{definition}[theorem]{\textbf{Definition}}
\newtheorem{remark}[theorem]{\textbf{Remark}}

\allowdisplaybreaks
\newcommand{\beq}{\begin{equation}}
\newcommand{\eeq}{\end{equation}}

\def \rz { {\mathbb R}}

\def \rz {{\mathbb R}}

\newcommand {\ar}{\rightarrow}
\newcommand {\pa}{\partial}
\newcommand{\qed}{$\Box$} %% Figure d\'{e}j\`{a} dans le package amsthm
\numberwithin{equation}{section}
\title{Remarks on the boundary set of  spectral equipartitions}
\author{P. B\'erard \\ Institut Fourier, Universit\'{e} Grenoble 1 and
CNRS, B.P.74,
\\ F 38 402 Saint Martin d'H\`{e}res Cedex, France.\\
and \\
B. Helffer \\
Laboratoire de Math\'ematiques, Univ Paris-Sud and CNRS,\\
F 91 405 Orsay Cedex, France.}
\date{\small{March 1, 2013 ~(to appear in Phil. Trans. A, Royal Society)}}

\begin{document}
\maketitle

\section{Introduction}
%% Ancienne introduction
% Given a bounded  open set $\Omega$ in $\mathbb{R}^n$ (or in a
% Riemannian  manifold),  and a partition of $\Omega$ by $k$ open sets
% $\omega_j$, we consider the quantity $\max_j \lambda(\omega_j)$,
% where $\lambda(\omega_j)$ is the ground state energy of the
% Dirichlet realization of the Laplacian in $\omega_j$.
%
% We denote by $\mathfrak{L}_k(\Omega)$ the infimum of $\max_j
% \lambda(\omega_j)$ over all $k$-partitions. A minimal $k$-partition
% is a partition which realizes the infimum. Although the analysis of
% minimal $k$-partitions is rather standard when $k=2$ (we find the
% nodal domains of  a second eigenfunction), the analysis for higher
% values of $k$ becomes non trivial and quite interesting. Minimal
% partitions are in particular spectral equipartitions, \emph{i.e.}
% the ground state energies $\lambda(\omega_j)$ are all equal.
%
% The purpose of this paper is to revisit various properties of nodal
% sets, and to explore if they are also true for minimal partitions,
% or more generally for spectral equipartitions. In
% Section~\ref{s-lb}, we prove a lower bound for the length of the
% boundary set of a partition in the $2$-dimensional situation. In
% Section~\ref{s-kl}, we consider estimates involving the cardinality
% of the partition.
%% Nouvelle introduction

Let $\Omega$ be a bounded open set in $\mathbb{R}^n$ (or in a
Riemannian manifold), and $k$ a positive integer $k$.

We denote by $\mathcal{D}_k(\Omega)$ the set of $k$-tuples
$\mathcal{D} = \{\omega_1, \ldots , \omega_k\}$ of open, pairwise
disjoint subsets of $\Omega$. Denote by $\phi$ a positive function
on open subsets of $\Omega$, and assume that it satisfies the
following two basic assumptions (i) if $\omega_1 \subset \omega_2$,
then $\phi(\omega_2) \le \phi(\omega_1)$, with strict inequality
provided that $\omega_2\setminus \omega_1$ is not too small; (ii)
$\phi(\omega)$ tends to infinity if $\omega$ skrinks to a point.
Introduce the function $\Phi$ defined on $\mathcal{D}_k(\Omega)$ by
$\Phi(\mathcal{D}) = \max_{1\le j \le k} \phi(\omega_k)$, if
$\mathcal{D} = \{\omega_1, \ldots , \omega_k\}$. Minimizing the
function $\Phi$ over $\mathcal{D}_k(\Omega)$ appears naturally when
one considers $k$ non-mixing populations (represented by the sets
$\omega_j$) competing for space (with the constraint represented by
the function $\phi$). As possible examples of constraint functions,
we mention the functions $\omega \to \mathrm{Vol}(\omega)^{-1}$ and
$\omega \to \lambda_1(\omega)$, where $\lambda_1(\omega)$ is the
ground state energy of the Dirichlet realization of the Laplacian in
$\omega$. It is easy to see that if there exists a minimal $k$-tuple
$\mathcal{D} = \{\omega_1, \ldots, \omega_k\}$, \emph{i.e.} one
which realizes the infimum of $\Phi$ over $\mathcal{D}_k(\Omega)$,
then $\mathcal{D}$ must actually be a partition of $\Omega$ and,
more precisely, an equipartition, which means that $\mathcal{D}$
also satisfies $\phi(\omega_1) = \cdots = \phi(\omega_k)$. For more
details on this subject, we refer to the papers \cite{BBO, CL1,
CTV0, CTV2, CTV:2005, Ha, Hel}.

In this paper, we consider the constraint function $D \to
\lambda_1(D)$. In this framework, it can be shown that minimal
partitions do exist, and that they have regular representatives (see
Properties~\ref{PP-1}). Although the analysis of minimal
$k$-partitions is rather standard when $k=2$ (we find the nodal
domains of  a second eigenfunction), the analysis for higher values
of $k$ becomes non trivial and quite interesting. If $u$ is an
eigenfunction of the Dirichlet realization of the Laplacian in
$\Omega$, the nodal domains of $u$ form an equipartition of $\Omega$
with boundaries the nodal set of $u$.

The purpose of this paper is to revisit various properties of nodal
sets, and to explore if they are also true for minimal spectral
partitions, or more generally for spectral equipartitions (here the
word \emph{spectral} refers to our choice of constraint function).

In Section~\ref{s-ini}, we review some basic notions related to our
choice of constraint function. In Section~\ref{s-lb}, we prove a
lower bound for the length of the boundary set of a partition in the
$2$-dimensional situation. In Section~\ref{s-kl}, we consider
estimates involving the cardinality of the partition.

\section{Definitions and notations}\label{s-ini}

\subsection{Spectral theory}

Let $\Omega$ be a bounded domain in $\mathbb{R}^2$, or a compact
Riemannian surface, possibly with boundary $\partial \Omega$, which
we assume to be piecewise $C^1$.  Let $H(\Omega)$ be the realization
of the Laplacian, or of the Laplace-Beltrami operator, $- \Delta$ in
$\Omega$, with Dirichlet boundary condition ($u|_{\partial \Omega
}=0)$. Let $\{\lambda_j(\Omega)\}_{j\ge 1}$ be the increasing
sequence of the eigenvalues of $H(\Omega)$, counted with
multiplicity. The eigenspace associated with $\lambda_k$ is denoted
by $E(\lambda_k)$. \medskip

A groundstate $u \in E(\lambda_1)$ does not vanish in $\Omega$ and
can be chosen to be positive. On the contrary, any non-zero
eigenfunction $u \in E(\lambda_k), k\ge 2$, changes sign in
$\Omega$, and hence has a nonempty zero set or \emph{nodal set},
\begin{equation}
N(u)=\overline{\{x\in \Omega\:\big|\: u(x)=0\}}.
\end{equation}
The connected components of $\Omega\setminus N(u)$ are called the
\emph{nodal domains} of $u$. The number of nodal domains of $u$ is
denoted by $\mu(u)$.\medskip

\noindent Courant's nodal domain theorem says:
\begin{theorem}[Courant]
Let $k\geq 1$, and let $E(\lambda_k)$ be the eigenspace of
$H(\Omega)$ associated with the eigenvalue $\lambda_k$. Then,
$\forall u\in E(\lambda_k)\setminus\{0\}\,,\; \mu (u)\le k\,.$
\end{theorem}

\noindent Except in dimension $1$, the inequality is strict in general. More
precisely, we have:
\begin{theorem}[Pleijel]
Let $\Omega$ be a bounded domain in $\mathbb{R}^2$. There exists a
constant $k_0$ depending on $\Omega$, such that if $k\geq k_0$, then
$$ \mu(u) < k\,,\,   \forall u\in  E(\lambda_k)\setminus \{0\}\,.
$$
\end{theorem}
Both theorems are proved in  \cite{Pleijel:1956}. The main points in
the proof of Pleijel's Theorem are the Faber-Krahn inequality and
the Weyl asymptotic law. Faber-Krahn's inequality states that, for
any bounded domain $\omega$ in $\mathbb{R}^2$,
\begin{equation}\label{E-fki}
\lambda_1(\omega) \geq \frac{\pi {\bf j}^2}{A(\omega)}\,,
\end{equation}
where $A(\omega)$ is the area of $\omega$, and ${\bf j}$ is the
least positive zero of the Bessel function of order $0$ (${\bf j}
\sim 2.4048$). Weyl's asymptotic law for the eigenvalues of
$H(\omega)$ states that
\begin{equation}\label{E-wal}
\lim_{k\to \infty} \frac{\lambda_k(\omega)}{k} =
\frac{4\pi}{A(\omega)}\,.
\end{equation}

\noindent Let $\bar \mu(k)$ be the maximum value of $\mu(u)$ when $u \in
E(\lambda_k)\setminus\{0\}$. Combining the results of Faber-Krahn
and Weyl, we obtain,
\begin{equation}\label{E-pleijel}
\limsup_{k\rightarrow +\infty} \frac{\bar \mu(k)}{k} \leq 4/{\bf j} ^2 <1 \,.
\end{equation}

 \begin{remark} Pleijel's
Theorem extends to bounded domains in $\mathbb{R}^n$, and more
generally to compact $n$-manifolds with boundary, with a constant
$\gamma(n) <1$ replacing $ 4/{\bf j} ^2$ in the right-hand side of
(\ref{E-pleijel}) (Peetre \cite{Pe}, B\'{e}rard-Meyer \cite{BeMe}). It
is also interesting to note that this constant is independent of the
geometry.
\end{remark}

\begin{remark} It follows from Pleijel's Theorem that the
equality $\bar{\mu}(k)=k$ can only occur for finitely many values of
$k$. The analysis of the equality  case  is very interesting. We
refer to \cite{HHOT} for more details.\end{remark}

\begin{remark} In dimension $1$, counting the nodal domains of
an eigenfunction of a Dirichlet Sturm-Liouville problem in some
interval $[a,b]$ is the same as counting the number of zeroes of the
eigenfunction. An analog in dimension $2$ is to consider the length
of the nodal set of eigenfunctions instead of the number of their
nodal domains. We shall come back to this question in
Section~\ref{s-lb}.
\end{remark}

\subsection{Partitions}\label{ss-p}
For this section, we refer to \cite{HHOT}. Let $k$ be a positive
integer. A (weak)  $k$-\emph{partition} of the open bounded set
$\Omega$ is\footnote{Note that we start from a very weak notion of
partition. We refer to \cite{HHOT} for a more precise definition of
classes of $k$-partitions, and for the notion of regular
representatives.} a family $\mathcal{D} = \{D_j\}_{j=1}^{k}$ of
pairwise disjoint sets such that $\cup_{j=1}^{k} D_j \subset
\Omega$. We denote by $\mathfrak{D}_k = \mathfrak{D}_k(\Omega)$ the
set of $k$-partitions such that the domains $D_j$ are open and
connected.\medskip

Given $\mathcal{D} \in \mathfrak{D}_k$, we define the \emph{energy}
$\Lambda(\mathcal{D})$ of the partition as,
\begin{equation}
\Lambda(\mathcal{D}) = \max_j \lambda(D_j),
\end{equation}
where $\lambda(D_j)$ is the groundstate energy of $H(D_j)$. We now
define the number $\mathfrak{L}_k(\Omega)$ as,
\begin{equation}
\mathfrak{L}_k(\Omega) = \inf_{\mathcal{D} \in \mathfrak{D}_k}
\Lambda(\mathcal{D})\,.
\end{equation}

\noindent A partition $\mathcal{D}$ is called \emph{minimal} if {$\Lambda(\mathcal{D}) =
\mathfrak{L}_k(\Omega).$} \medskip

\noindent \textbf{Example}. The nodal domains of an eigenfunction $u
\in E(\lambda)\setminus\{0\}$ of $H(\Omega)$ form a
$\mu(u)$-partition of $\Omega$ denoted by $\mathcal D(u)$. Such a
partition is called a \emph{nodal partition}.\medskip

\noindent It turns out that $\mathfrak{L}_2(\Omega) = \lambda_2(\Omega)$, and
that minimal $2$-partitions are nodal partitions. The situation when
$k \ge 3$ is more complicated, and more interesting,
\cite{HHOT}.\medskip

A partition $\mathcal{D} = \{D_j\}_{j=1}^{k} \in
\mathfrak{D}_k(\Omega)$ is called \emph{strong} if,
\begin{equation}
\mathrm{Int}(\overline{\cup_j D_j} ) = \Omega
\,.
\end{equation}

The \emph{boundary set} $N(\mathcal{D})$ of a strong partition
$\mathcal{D}_k = \{D_j\}_{j=1}^{k} \in \mathfrak{D}(\Omega)$ is the
closed set,
\begin{equation}
N(\mathcal{D}) = \overline{\cup_j(\partial D_j \cap \Omega)}\,.
\end{equation}

The set $\mathcal{R}_k(\Omega)$ of \emph{regular $k$-partitions} is
the subset of strong $k$-partitions in $\mathfrak{D}_k(\Omega)$
whose boundary set $N = N(\mathcal{D})$ satisfies the following
properties.
\begin{itemize}
    \item[(i)] The set $N$ is locally a regular curve in $\Omega$,
    except possibly at finitely many points $\{y_i\} \in N \cap
    \Omega$, in the neighborhood of which $N$ is the union of
    $\nu(y_i)$ smooth semi-arcs at $y_i, \nu(y_i) \ge 3$.
    \item[(ii)] The set $N\cap \partial \Omega$ consists of finitely
    many points $\{z_j\}$. Near the point $z_j$, the set $N$ is the
    union of $\rho(z_j)\ge 1$ semi-arcs hitting $\partial \Omega$ at $z_j$.
    \item[(iii)] The set $N$ has the equal angle property. More precisely,
    at any \emph{interior singular point} $y_i$, the semi-arcs meet with equal
    angles; at any \emph{boundary singular point} $z_j$, the semi-arcs
    form equal angles together with the boundary $\partial \Omega$.
\end{itemize}\medskip

\noindent \textbf{Example}. A nodal partition $\mathcal{D}(u)$
provides an example of a regular partition, and the boundary set
$N\big(\mathcal{D}(u)\big)$ coincides with the nodal set $N(u)$.
Note that for a regular partition, the number $\nu(y_i)$ of
semi-arcs at an interior singular point may be odd, whereas it is
always even for a nodal partition.\medskip

Let us now introduce:

\begin{definition}
We call \emph{spectral equipartition} a strong $k$-partition
$\mathcal{D} =\{D_i\}$ such that $\lambda(D_i) =
\Lambda(\mathcal{D})\,$, for $ i=1, \dots, k$. The number $
\Lambda(\mathcal D)$ is called the \emph{energy} of the
equipartition.
\end{definition}

\noindent \textbf{Example}. Nodal partitions provide examples of
spectral equipartitions.

\begin{properties}\label{PP-1}
Given an open bounded set $\Omega$,
\begin{enumerate}
    \item[(i)] Minimal $k$-partitions exist (\cite{CTV0, CTV2,
    CTV:2005}).
    \item[(ii)] Any minimal $k$-partition has a representative (modulo
    sets of capacity $0$) which is a regular spectral equipartition
    (\cite{HHOT}).
\end{enumerate}
\end{properties}

\subsection{Euler formula}
Let $\Omega \subset M$ be a bounded domain with piecewise smooth
boundary  $\partial \Omega$. Let $N \subset \overline{\Omega}$ be a
regular closed set (in the sense of Section~\ref{ss-p}, properties
(i)-(iii)) such that the family $\mathcal{D} = \{D_1, \ldots, D_k\}$
of connected components of $\Omega\!\setminus\! N$ is a regular,
strong partition of $\Omega$. Recall that for a singular point $y
\in N \cap \Omega$, $\nu(y)$ is the number of semi-arcs at $y$, and
that for a singular point $z \in N \cap \partial \Omega$, $\rho(z)$
is the number of semi-arcs at $z$, not counting the two arcs
contained in $\partial \Omega$. Let $\mathcal{S}(\mathcal{D})$
denote the set of singular points of $N(\mathcal{D})$, both interior
or boundary points, if any. We define the \emph{index} of a point $x
\in \mathcal{S}(\mathcal{D})$ to be,
\begin{equation}
\iota (x) = \left\{
\begin{array}{ll}
\nu(x) - 2\,, & \text{~if~} x \text{ ~is an interior singular point},\\
\rho(x)\,,  & \text{~if~} x \text{ ~is a boundary singular point}.
\end{array}
\right.
\end{equation}
We introduce the number $\sigma(\mathcal{D})$ to be,
\begin{equation}
\sigma(\mathcal{D}) = \sum_{x \in \mathcal{S}(\mathcal{D})}
\iota(x)\,.
\end{equation}
For a regular strong $k$-partition $\mathcal{D} = \{D_j\}_{j=1}^k$
of $\Omega$, we have Euler's formula,
\begin{equation}\label{E-euler}
\chi(\Omega) + \frac{1}{2} \sigma(\mathcal{D}) = \sum_{j=1}^k
\chi(D_j)\,.
\end{equation}

We refer to  \cite{HOMN} for a combinatorial proof of this formula in the
case of an open set of $\mathbb R^2$. One can give a Riemannian
proof using the global Gauss-Bonnet theorem. For a domain $D$ with
piecewise smooth boundary $\partial D$ consisting of piecewise $C^1$
simple closed curves $\{C_i\}_{i=1}^{n}$, with corners
$\{p_{i,j}\}$ ($i=1, \dots, n$ and  $j=1,\dots, m_i$) and corresponding interior
angles $\theta_{i,j}\,$, we have
\begin{equation}
2\pi \chi(D) = \int_D K  + \sum_{i=1}^n \beta(C_i)\,,
\end{equation}
where
$$\beta(C_i) = \int_{C_i}\langle k,\nu_{D}\rangle  + \sum_{j=1}^{m_i}
(\pi - \theta_{i,j})\,.
$$
In this formula, $k$ is the geodesic curvature vector of the regular
part of the curve $C_i$, and $\nu_{D}$ is the unit normal to $C_i$
pointing inside $D$.\medskip

\noindent To prove (\ref{E-euler}), it suffices to sum up the
Gauss-Bonnet formulas relative to each domain $D_j$ and to take into
account the following facts:
\begin{itemize}
    \item the integrals of the Gaussian curvature over the $D_j$'s
    add up to the integral of the Gaussian curvature over
    $\Omega$,
    \item cancellations occur when adding the integrals of the
    geodesic curvature over the curves bounding two adjacent $D_j$ (the unit
    normal vectors point in opposite directions), while they add up to give the
    integral of the geodesic curvature over the boundary of $\Omega$,
    \item there are contributions coming from the angles associated
    with the singular points of $N$ and, when summed up, these contributions yield
    the second term in the left-hand side of (\ref{E-euler}).
\end{itemize}

\noindent Note that the proof of (\ref{E-euler}) does not use the
fact that the semi-arcs meet at the singular points of $N$ with
equal angles.

\section{Lower bounds for the length of the boundary set of a
regular spectral equipartition}\label{s-lb}

\subsection{Introduction}\label{ss-into}

Let $\mathcal{D} = \{D_1, \ldots , D_k\}$ be a regular
spectral equipartition with energy $\Lambda =
 \Lambda(\mathcal{D})$. The boundary set $N(\mathcal{D})$ of the
partition consists of singular points $\{y_i\}_{i=1}^a$ inside
$\Omega$, of singular points $\{z_i\}_{i=1}^b$ on $\partial \Omega$,
of $C^1$ arcs $\{\gamma_i\}_{i=1}^c$ which bound two adjacent
domains of the partition, and of arcs $\{\delta_i\}_{i=1}^d$
contained in $\partial \Omega$. We define the length of the boundary
set $N(\mathcal{D})$ by the formula,
\begin{equation}\label{E-lbs}
P(\mathcal{D}) := \sum_{i=1}^c \ell(\gamma_i) + \frac{1}{2}
\ell(\partial \Omega),
\end{equation}
where $\ell$ denotes the length of the curves. Note that
$\ell(\partial \Omega) = \sum_{i=1}^d \ell(\delta_i)$
and that $\sum_{i=1}^k \ell (\partial D_i) = 2
P(\mathcal{D})$.\medskip

In this section, we investigate lower bounds for $P(\mathcal{D})$ in
terms of the energy $\Lambda(\mathcal{D})$ and the area $A(\Omega)$.
\medskip

\noindent As a matter of fact, we show that the methods introduced in
\cite{BrGr, Br, Sa1} apply to regular spectral equipartitions, and
hence to minimal partitions. We provide three estimates.
\begin{enumerate}
    \item The first estimate holds for plane domains, and follows the
    method of \cite{BrGr}.
    \item The second estimate applies to a compact Riemannian
    surface (with or without boundary), and follows the method of
    \cite{Sa1}.
    \item The third estimate is a local estimate based on the method
    of \cite{Br}.
\end{enumerate}

\noindent Let $\mathcal{D} = \{D_i\}_{i=1}^{k}$ be a regular spectral
equipartition with energy $\Lambda =
\Lambda(\mathcal{D})$. Let $R(D_i)$ be the inner radius of the set
$D_i$. Recall that ${\bf j}$ denotes the least positive zero of the
Bessel function of order $0$.

\subsection{The method of Br\"uning-Gromes}

In this section, $\Omega$ is a bounded domain in $\mathbb{R}^2$,
with piecewise $C^1$ boundary. We only sketch the method which
relies on three inequalities.

\begin{enumerate}
    \item The monotonicity of eigenvalues and the characterization
    of the ground state imply that
    \begin{equation}\label{bg-1}
    \forall i, 1 \le i \le k, ~~~R(D_i) \le \frac{{\bf j}}{\sqrt{\Lambda}}\,.
    \end{equation}
    \item The Faber-Krahn inequality and the isoperimetric
    inequality imply that
    \begin{equation}\label{bg-2}
    \forall i, 1 \le i \le k, ~~~\frac{2\pi{\bf  j}}{\sqrt{\Lambda}} \le
    \ell (\partial D_i)\,.
    \end{equation}
    \item The generalized F\'{e}jes-Toth isoperimetric inequality
    (\cite{BrGr}, Hilfssatz 2) asserts that, for $1 \le i \le k$,
    \begin{equation}\label{bg-3}
     A(D_i) \le R(D_i) \ell(\partial D_i) -
    \chi(D_i) \pi R^2(D_i)\,.
\end{equation}
\end{enumerate}
Using that $\chi(D_i)\leq 1$, we immediately see that the function
$r \to r \ell(\partial D_i) -\chi (D_i)  \pi r^2$ is non-decreasing
for $2\pi \, r \le \ell(\partial D_i)$. Using inequalities
(\ref{bg-1}) and (\ref{bg-2}), it follows that one can substitute
$\frac{j}{\sqrt{\Lambda}}$ to $R(D_i)$ in \eqref{bg-3} and obtain,
\begin{equation}\label{bg-4}
 A(D_i) \le \frac{\mathbf{j}}{\sqrt{\Lambda}} \ell(\partial
D_i) - \chi(D_i) \pi \big( \frac{\mathbf{j}}{\sqrt{\Lambda}}
\big)^2, \text{ ~for~ } 1 \le i \le k.
\end{equation}

\noindent Summing up the inequalities \eqref{bg-4}, for $1\le i \le k$, we
obtain
$$
A(\Omega ) \le \frac{{\bf j}}{\sqrt{\Lambda}} \sum_i \ell (\partial D_i)
- \sum_i \chi(D_i) \pi \frac{{\bf j}^2}{\Lambda}\,.
$$
Using Euler's formula (\ref{E-euler}), we conclude that
\begin{equation*}
A(\Omega )  \le \frac{2{\bf j}}{\sqrt{\Lambda}} P(\mathcal{D}) -
\big[\chi(\Omega ) + \frac{1}{2} \sigma(\mathcal{D})\big]\pi
\frac{{\bf j}^2}{\Lambda}\,.
\end{equation*}

\noindent We have proved,

\begin{proposition}\label{P-bg}
Let $\Omega$ be a bounded open set in $\rz^2$, and let $\mathcal{D}$
be a regular spectral equipartition of $\Omega$. The length
$P(\mathcal{D})$ of the boundary set of $\mathcal{D}$ is bounded
from below in terms of the energy $\Lambda(\mathcal{D})$. More
precisely,
\begin{equation}\label{bg}
\frac{A(\Omega ) }{2{\bf j} } \sqrt{\Lambda(\mathcal{D})}  +
\frac{\pi{\bf j}}{2\sqrt{\Lambda(\mathcal{D})}} \big[\chi(\Omega ) +
\frac{1}{2} \sigma(\mathcal{D})\big] \le P(\mathcal{D})\,.
\end{equation}
\end{proposition}

Note that (\ref{bg}) is actually slightly better than the estimate
in \cite{BrGr} which does not take into account the term
$\sigma(\mathcal{D})$ when $\mathcal{D}$ is the nodal partition for
an eigenfunction $u$ associated with the eigenvalue $\Lambda$. This
fact is suggested in \cite{Sa1}.

\subsection{The method of Savo}

In this section, we follow the method of Savo \cite{Sa1}, and keep
the same notations and assumptions. We sketch the proof in the case
with boundary as it is not detailed in \cite{Sa1}. Here, $\Omega$ is
a compact Riemannian surface with boundary.  We denote the
Laplace-Beltrami operator by $\Delta$ and the Gaussian curvature by
$K$. We write $K = K_+ - K_-$ (the negative and positive parts of
the curvature). \medskip

\noindent We assume that numbers $\alpha \geq 0$ and $D$ are
given such that:
$$K \ge - \alpha^2 \text{ ~and~ } \delta(\Omega) \leq D\,\,,$$
where $\delta(\Omega)$ is the diameter of $\Omega$. Finally, we
define the numbers $$B(\Omega) = \int_{\Omega} K_+ - 2\pi
\chi(\Omega)\,,$$ and $$C(\alpha,D) = \sqrt{\pi^2 + \frac{1}{4}
\alpha^2 D^2}\,.
$$
We recall the following results from \cite{Sa1}.

\begin{lemma}\label{L-A}
Let $\Omega$ be a compact Riemannian surface with piecewise $C^1$ boundary. Then
\begin{equation}\label{savo1}
\frac{2}{\pi} A(\Omega) \sqrt{\lambda(\Omega)} \le \ell(\partial \Omega) + R(\Omega)
\max\{B(\Omega),0\}\,,
\end{equation}
where $\lambda(\Omega)$ is the ground state energy of the Dirichlet
realization fo the Laplacian in $\Omega$, and $R(\Omega)$ the inner
radius of $\Omega$.
\end{lemma}

\noindent This is Proposition~3 in \cite{Sa1} (p.~137). Note that when $M$ is flat
and $\Omega$ is simply
or doubly connected, we recover Polya's inequality \cite{Pol}  which reads:
\begin{equation}\label{Polya}
\frac{2}{\pi} A(\Omega) \sqrt{\lambda(\Omega)} \le \ell(\partial \Omega)\,.
\end{equation}

\begin{lemma}\label{L-B}
Let $\Omega$ be a compact Riemannian surface with piecewise $C^1$ boundary. Then,
$$R(\Omega) \sqrt{\lambda(\Omega )} \le \min \Big\{C(\alpha,D), \sqrt{\pi^2 +
\frac{\alpha^2 C^2(\alpha ,D)}{4 \lambda(\Omega)}} \Big\} =:
\psi(\alpha ; D ; \lambda(\Omega)).$$
\end{lemma}

\noindent This is Lemma~10 in \cite{Sa1} (p.~141) using $\lambda(\Omega)$
instead of $\lambda$.

\begin{lemma}\label{L-C}
Let $\Omega$ be a compact Riemannian surface with piecewise $C^1$
boundary. Assume that $B(\Omega) < 0$. Then,
$$2 |B(\Omega)| \le \lambda(\Omega) A(\Omega) \le \frac{\pi}{2} \sqrt{\lambda(\Omega)}
\ell(\partial \Omega)\,.$$ \end{lemma}

\noindent This is Lemma~11 in \cite{Sa1} (p. 141), which relies on Dong's
paper \cite{Dong}. Note that the second inequality follows from
Lemma~\ref{L-A}.
\medskip

Let us now proceed with the lower estimate of
$P(\mathcal{D})$ when $\Omega$ is a Riemannian surface with
boundary.

\begin{proposition}\label{P-lower}
Let $\Omega$ be a compact Riemannian surface with piecewise $C^1$
boundary. The length $P(\mathcal{D})$ of the boundary set of a
regular spectral equipartition $\mathcal{D}$, with
energy $\Lambda(\mathcal{D})=\Lambda$, satisfies the
inequality
\begin{equation}\label{E-lower}
P(\mathcal{D}) \ge \frac{4A(\Omega) \sqrt{\Lambda}}{4\pi +
\pi^2\psi(\alpha,D;\Lambda)} - \frac{2\pi \psi(\alpha,D;\Lambda)}{
\sqrt{\Lambda}\big( 4\pi + \pi^2 \psi(\alpha,D;\Lambda) \big)}\big(
B(\Omega) - \pi \sigma(\mathcal{D})\big).
\end{equation}
\end{proposition}

\noindent \textbf{Proof}. The proof follows the ideas in \cite{Sa1} closely.
Since Savo does not provide all the details for the case with
boundary, we provide them here. Lemma~\ref{L-A} applied to each
$D_j$ gives,
$$\frac{2}{\pi} A(D_j) \sqrt{\lambda(D_j)} \le \ell(\partial D_j) + R(D_j)
\max\{B(D_j),0\}.$$

\noindent Since $\lambda(D_j)=\Lambda$ for all $j$, summing up in $j$, we find
that

\begin{equation}\label{E-a}
\frac{2}{\pi} A(\Omega) \sqrt{\Lambda} \le 2 P(\mathcal{D}) +
\sum_{j=1}^k R(D_j) \max\{B(D_j),0\}.
\end{equation}

\noindent Call $T$ the second term in the right-hand side of the preceding
inequality and define the sets,
$$J_+ := \{j ~|~ 1 \le j \le k, ~B(D_j) > 0\}, ~~ J_- := \{j ~|~ 1 \le j
\le k, ~B(D_j) \le 0\}.$$

By Lemma~\ref{L-B}, we have
\begin{equation} T = \sum_{j \in J_+} R(D_j) B(D_j) \le \frac{\psi(\alpha , D ;
\Lambda)}{\sqrt{\Lambda}}
\sum_{j \in J_+} B(D_j).
\end{equation}

\noindent Using the definition of $B(D_j)$, we find that
$$
\begin{array}{ll}
\sum_{j=1}^k B(D_j) & = \int_{\Omega} K_+  - 2\pi \sum_{j=1}^k
\chi(D_j) \\[6pt]
    & = B(\Omega) + 2\pi \chi(\Omega)  - 2\pi \sum_{j=1}^k
\chi(D_j) \\
\end{array}
$$
and hence, using Euler's formula (\ref{E-euler}),
$$\sum_{j=1}^k B(D_j) = B(\Omega) -  \pi \sigma(\mathcal{D}).$$

\noindent On the other hand, we have
$$
\begin{array}{ll}
\sum_{j \in J_+} B(D_j) & = \sum_{j=1}^k B(D_j) - \sum_{j \in J_-} B(D_j) \\[6pt]
    & = B(\Omega) -  \pi \sigma(\mathcal{D}) + \sum_{j \in J_-}
    |B(D_j)|,
\end{array}
$$
and we can estimate the last term in the right-hand side using
Lemma~\ref{L-C}. Namely,
$$\sum_{j \in J_-}
    |B(D_j)| \le \frac{\pi}{4} \sqrt{\Lambda}\sum_{j \in J_-} \ell(\partial D_j) \le
    \frac{\pi}{2} \sqrt{\Lambda} P(\mathcal{D}).$$
Finally, we obtain the following estimate for $T$,

$$T \le \frac{\psi(\alpha , D ; \Lambda)}{\sqrt{\Lambda}}
\big\{ B(\Omega) -  \pi \sigma(\mathcal{D}) + \frac{\pi}{2}
\sqrt{\Lambda} P(\mathcal{D}) \big\}.$$

\noindent Using (\ref{E-a}), it follows that
\begin{equation}\label{E-b}
A(\Omega) \sqrt{\Lambda} \le \big\{ \pi +
\frac{\pi^2}{4}\psi(\alpha, D; \Lambda) \big\} P(\mathcal{D}) +
\frac{\pi \, \psi(\alpha , D ; \Lambda)}{2 \sqrt{\Lambda}} \big\{
B(\Omega) - \pi \sigma(\mathcal{D}) \big\}.
\end{equation}

\noindent This proves the proposition. \qed

\subsection{A loose local lower estimate for $P(\mathcal{D})$}

For simplicity, we now assume that $\Omega$ is a bounded domain in
$\mathbb{R}^2$, with piecewise $C^1$ boundary. We also assume that
we are given some point $x_0 \in \Omega$, some radius $R$ and some
positive number $\rho$, small with respect to $R$, such that
$B(x_0,R+\rho) \subset \Omega$. Note that the ball $B(x_0,R)$ could
be replaced by any regular domain.

\subsubsection{A local estimate \`{a} la Br\"{u}ning-Gromes : eigenvalues}

\begin{lemma}\label{L-A1}
Let $\lambda $ be an eigenvalue of $H(\Omega)$, and let $u \in
E(\lambda)$ be a non-zero eigenfunction associated with $\lambda$.
If $\lambda r^2 > {\bf j}^2$, then any disk $B(x,r) \subset \Omega$
contains at least a point of the nodal set $N(u)$.
\end{lemma}

\noindent This follows immediately from the monotonicity of the Dirichlet
eigenvalues with respect to domain inclusion.

\begin{lemma}\label{L-B1}
Let $\lambda$ be an eigenvalue of $H(\Omega)$,  assumed
to be large enough. Let $r>0$ be such that $0 < r \le \rho <
\frac{R}{10}$, and $\lambda r ^2 > 4 {\bf j}^2$. Then, there exists
a family of points $\{x_1, \ldots, x_{N}\}$ such that:
\begin{itemize}
    \item[(1)] For $1 \le j \le N$, $x_j \in N(u) \cap
    B(x_0, R-\frac{r}{2})$.
    \item[(2)] The balls $B(x_j, \frac{r}{2}), 1 \le j \le N$, are
    pairwise disjoint and contained in $B(x_0,R) \subset \Omega$.
    \item[(3)] We have the inclusion $B(x_0,R-r) \subset
    \cup_{j=1}^{N} B(x_j,2r)$.
    \item[(4)] The number $N$ satisfies, $r^2 N \ge 0.2\, R^2$.
\end{itemize}
\end{lemma}

\noindent \emph{Proof}. (a) Consider the ball $B(x_0,R-r)$ and take $y_1, y_2$
to be the end points of a diameter of the closed ball. Because $r
\le \rho < R/10$ and $r^2 \lambda  > 4 {\bf j}^2$, we have that
$B(y_i,\frac{r}{2}) \subset B(x_0,R-\frac{r}{2}) \subset \Omega$ and
$B(y_i,\frac{r}{2}) \cap N(u) \not = \emptyset$. Choose $x_i$ in
$B(y_i,\frac{r}{2})\cap N(u)$. Then, $x_i \in N(u) \cap
B(x_0,R-\frac{r}{2})$, $B(x_1,\frac{r}{2}) \cap B(x_2,\frac{r}{2}) =
\emptyset$ and $B(x_i,\frac{r}{2}) \subset B(x_0,R) \subset \Omega$.
\medskip

\noindent (b) Take a maximal element $\{x_1, \ldots , x_{N}\}$ (with respect
to inclusion) in the set
$$\mathcal{F} := \big\{ (x_1, \ldots, x_k) ~|~ x_i \in N(u) \cap B(x_0,R-\frac{r}{2}),
B(x_i, \frac{r}{2}) \text{ pairwise disjoint} \big\},$$ so that the
family $\{x_1, \ldots , x_{N}\}$ satisfies (1) and (2). \medskip

\noindent \emph{We claim that (3) holds}. Indeed, otherwise we could find $y
\in B(x_0,R-r)$ with $d(x_i,y) \ge 2r$, for $1 \le i \le N$. Because
$B(y,\frac{r}{2}) \cap N(u) \not = \emptyset$, we would find some $
z \in B(x_0,R-\frac{r}{2}) \cap N(u) \cap B(y,\frac{r}{2})$ such
that $B(z,\frac{r}{2}) \cap \big( \cup_{j=1}^{N} B(x_j,\frac{r}{2})
\big) = \emptyset$. This would contradict the maximality of the
family.
\medskip

\noindent (c) Assertion (3) implies that $\pi (R-r)^2 \le \sum_{j=1}^{N}
A(B(x_j, 2r)) = 4\pi N r^2$ and since $r \le \rho < R/10$, we get
$r^2 N \ge (0.9)^2 \frac{1}{4} R^2$. The Lemma is proved. \qed
\medskip

\noindent Recall that $N(u)$ consists of finitely many points and finitely
many $C^1$ arcs with finite length. \medskip

\begin{lemma}\label{L-C1}
Let $\{x_1, \ldots , x_{N}\}$ be a maximal family as given by
Lemma~\ref{L-B1}. Assume that $r^2 \lambda < 16\, {\bf j}^2$. Then
there exists no nodal curve $\gamma \subset N(u)$ which is simply
closed and contained in any of the balls $B(x_j, \frac{r}{4}), 1 \le
j \le N$.
\end{lemma}

\noindent \emph{Proof}. Indeed, otherwise, there would be a nodal domain
contained in one of the balls $B(x_j, \frac{r}{4})$ and hence we
would have $r^2 \, \lambda \ge 16\, {\bf j}^2$.\medskip

We can now prove the following local estimate.

\begin{proposition}\label{P-local1}
Let $\lambda $ be an eigenvalue of $H(\Omega)$,
assumed to be large enough. Let $u$ be a non-zero
eigenfunction associated with $\lambda$. Then, the length of the
nodal set $N(u)$ inside $B(x_0,R)$ is bounded from below by $10^{-2}
R^2 \sqrt{\lambda}$.
\end{proposition}

\noindent \textbf{Proof}. Choose $(r,\lambda)$ so that $4\,{\bf j}^2 < r^2
\lambda  < 16\, {\bf j}^2$, with $r \le \rho <R/10$. By
Lemma~\ref{L-B1}, the $N$ balls $B(x_j, \frac{r}{4})$ are pairwise
disjoint with center on $N(u)$. By Lemma~\ref{L-C1}, the length of
$N(u) \cap B(x_j, \frac{r}{4})$ is at least $\frac{r}{2}$. It
follows that
\begin{equation}\label{E-loose}
\ell\big(N(u) \cap B(x_0, R)\big) \ge \sum_j \ell\big( N(u)\cap
B(x_j,\frac{r}{4})\big) \ge N \frac{r}{2}\,,
\end{equation}
and the result follows in view of the estimates $r^2 N \ge 0.2\,  r^2$
and $r^2 \lambda < 16 \, {\bf j}^2$.\medskip

\begin{remark}
Proposition~\ref{P-local1} can be generalized to the case of a
compact Riemannian surface with or without boundary. In that case,
one needs to consider balls with radii less than the injectivity
radius of the surface, and replace the Faber-Krahn inequality by a
local Faber-Krahn inequality, using the fact that the metric can be
at small scale compared with a Euclidean metric (see \cite{Br} for
more details).
\end{remark}

\subsubsection{A local estimate \`{a} la Br\"{u}ning-Gromes : spectral equipartitions}

The above proof applies to a regular spectral equipartition of
energy $\Lambda$.  It is enough in the statements to replace the
nodal set $N(u)$ of $u$ by the boundary set $N\big(\mathcal{D}\big)$
of the partition $\mathcal{D}$. We just rewrite the first statement.

\begin{lemma}\label{L-A1bis}
Let $\Lambda$ be the energy of a regular spectral equipartition.  If
$\Lambda r^2 > {\bf j}^2$, then any disk $B(x,r) \subset \Omega$
contains at least one point of boundary set of the partition.
\end{lemma}

\noindent This follows immediately from the monotonicity of the Dirichlet
eigenvalues with respect to domain inclusion. \medskip

\section{Estimates involving the cardinality of the\\ partitions}\label{s-kl}

\noindent Let $\mathcal{D} = \{D_i\}_{i=1}^k$ be a partition of $\Omega$. We
call the number $k$ the \emph{cardinality} of the partition, and we
denote it by $k = \sharp(\mathcal{D})$.\medskip

\subsection{Estimates on the energy $\Lambda(\mathcal{D})$ and on
$\mathfrak L_k(\Omega)$}

\begin{proposition}\label{P-fk}~
\begin{enumerate}
\item[(i)] Let $\Omega$ be a bounded open subset of $\rz^2$. The energy
$\Lambda(\mathcal{D})$ of a partition $\mathcal{D}$ of $\Omega$
satisfies the inequality,
\begin{equation}\label{E-fke}
\Lambda(\mathcal{D}) \geq \frac{\pi {\bf j}^2}{A(\Omega)} \,
\sharp(\mathcal{D})\,.
\end{equation}
In particular, for any $k \ge 1$, we have the inequality,
\begin{equation}\label{E-fkl}
\mathfrak L_k (\Omega) \geq  \frac{\pi {\bf j}^2}{A(\Omega)} \, k
\,.
\end{equation}
\item[(ii)] Let $\Omega$ be a bounded open subset on a compact Riemannian surface.
Then,
\begin{equation}\label{E-glb}
A(\Omega) \liminf_{k \ar \infty} \frac{\mathcal{L}_k(\Omega)}{k} \ge
 \pi \mathbf{j}^2 \,.
\end{equation}
\end{enumerate}
\end{proposition}

\noindent \textbf{Proof}. Assertion (i) is an immediate consequence of the Faber-Krahn
inequality \eqref{E-fki}. To prove Assertion~(ii), we use the fact that on a
general compact surface $M$, we have the following asymptotic isoperimetric
and Faber-Krahn inequalities (which actually hold in arbitrary dimension).

\begin{lemma}[\cite{BeMe}, Lemma~II.15, p. 528]
Let $(M,g)$ be a compact Riemannian surface. For any $\epsilon
> 0$, there exists a positive number $a(M,g,\epsilon)$ such that for
any regular domain $\omega \subset M$ with area $A(\omega)$
less than or equal to $a(M,g,\epsilon)$,
$$
\left\{%
\begin{array}{lll}
\ell(\partial \omega) & \ge & (1-\epsilon) \ell(\partial
\omega^{*}),\\
\lambda(\omega) & \ge & (1-\epsilon)^2 \frac{\pi
\mathbf{j}^2}{A(\omega)},
\end{array}%
\right.
$$
where $\omega^{*}$ is a Euclidean disk of area $A(\omega)$.
\end{lemma}

\noindent Let $\mathcal{D} =\{D_i\}$ be a partition of $\Omega$. Let
$$
J_{\epsilon} = \big\{i \in \{1, \ldots, k\} ~|~ A(D_i)
> a(M,g,\epsilon) \big\}.
$$
The number of elements of this set is bounded by
\begin{equation}\label{E-jeps}
\sharp(J_{\epsilon}) \le \frac{A(\Omega)}{a(M,g,\epsilon)}.
\end{equation}

\noindent For any $i \not \in J_{\epsilon}$, we can write,
$$
\lambda(D_i) \ge (1-\epsilon)^2 \frac{\pi \mathbf{j}^2}{A(D_i)}
$$
and hence, provided that $\sharp(\mathcal{D})$ is large enough,
$$
\Lambda(\mathcal{D}) A(\Omega) \ge (1-\epsilon)^2 (\sharp(\mathcal{D}) -
\frac{A(\Omega)}{a(M,g,\epsilon)}) \pi \mathbf{j}^2.
$$

\noindent As a consequence, we obtain that
$$
A(\Omega) \liminf_{k\ar  \infty} \frac{\mathcal{L}_k(\Omega)}{k} \ge
(1-\epsilon)^2 \pi \mathbf{j}^2.
$$

\noindent We can now let $\epsilon$ tend to zero to get the estimate
\eqref{E-glb}. \qed\medskip

\noindent \textbf{Remarks}.\\
(1) We point out that the lower bounds in the Proposition only depend on
the area of $\Omega$, not on its geometry.\\
(2) Similar inequalities on $\mathcal{L}_k(\Omega)$ can also be deduced
from \cite{Pe}, when $\Omega$ is a bounded domain in a simply-connected
surface $M$ with Gaussian curvature $K$, such that $\Omega \subset \Omega_0$, a
simply-connected domain satisfying $A(\Omega_0) \sup_{\Omega_0}K^+
\le \pi$. Let us mention two particular cases.\\
\noindent ~(a) If $M$ is a simply-connected surface with non-positive
curvature, then according to \cite{Pe}, $\lambda(\Omega)
A(\Omega) \ge \pi \mathbf{j}^2$ for any bounded domain $\Omega$
and we conclude that $$A(\Omega) \frac{\mathcal{L}_k(\Omega)}{k} \ge
\pi \mathbf{j}^2 $$ for all $k \ge 1$, as in the Euclidean case.\\
\noindent ~(b) If $M$ is the standard sphere, then according to \cite{Pe},
$$
\lambda(D) A(D) \ge \pi \mathbf{j}^2 \Big(1 - \frac{A(D)}{4\pi}\Big)
$$
for any domain $D$, and one can conclude that, for any domain $\Omega$,
$$
A(\Omega)
    \frac{\mathcal{L}_k(\Omega)}{k} \ge \pi \mathbf{j}^2 -
    \frac{\mathbf{j}^2 }{4k} A(\Omega)\,,
$$
for all $k \ge 1$.\medskip

\noindent \begin{remark}\label{R-fg1}
Given a $k$-partition $\mathcal{D} = \{D_i\}_{i=1}^k$ of $\Omega$,
one can also introduce the energy $\Lambda_1(\mathcal{D}) :=\frac{1}
{k}\,\sum_{i=1}^{k} \lambda (D_i)$, and define the number $\mathfrak
L_{k,1}(\Omega)$ by taking the infimum of the $\Lambda_1$-energy
over all $k$-partitions. An easy convexity argument shows that the
above inequalities \eqref{E-fke}, \eqref{E-fkl} and \eqref{E-glb}
hold with $\Lambda(\mathcal{D})$ and $\mathfrak{L}_k (\Omega)$ replaced by
$\Lambda_1(\mathcal{D})$ and $\mathfrak{L}_{k,1}(\Omega)$
respectively. For example, we have the inequalities,
\begin{equation*}
\Lambda(\mathcal{D}) \ge \Lambda_1(\mathcal{D}) \geq \frac{\pi {\bf
j}^2}{A(\Omega)} \, \sharp(\mathcal{D})\,, \text{ ~and,}
\end{equation*}
\begin{equation*}
\mathfrak L_k (\Omega) \ge \mathfrak L_{k,1}(\Omega) \geq  \frac{\pi
{\bf j}^2}{A(\Omega)} \, k \,.
\end{equation*}
\end{remark}\medskip

For bounded domains in $\rz^2$, one can also give an asymptotic
upper bound for $\mathfrak{L}_k$. More precisely,

\begin{property}\label{PP-2}
For any regular bounded open sutset of $\rz^2$,
\begin{equation}\label{E-asym-l}
\limsup_{k\rightarrow +\infty} \frac{\mathfrak L_k(\Omega)}{k}\ \leq
\frac{\lambda(\mathrm{Hexa}_1)}{A(\Omega)}\,,
\end{equation}
where $\mathrm{Hexa}_1$ is the regular hexagon in $\mathbb{R}^2$,
with area $1$.
\end{property}

This can be seen by considering the hexagonal tiling in the plane,
with hexagons of area $a$ and the partition of $\Omega$ given by
taking the union of the hexagons contained in $\Omega$, whose number
is asymptotically $\frac{A(\Omega)}{a}$ when $a$ tends to
zero.\medskip

The inequalities \eqref{E-asym-l} and \eqref{E-glb} motivate the following
two conjectures\footnote{The second author was informed of these conjectures
by M. Van den Berg.} for bounded domains in $\rz^2$.  They were proposed
and analyzed in the recent years (see \cite{BHV,BBO,CL1,HHOT,Hel}). The first
one is that
\begin{conjecture}\label{ConjAs1}
The limit of ${ {\mathfrak L_k(\Omega)}/{k}}$ as ${ k\rightarrow
+\infty}$  exists.
\end{conjecture}

\noindent The second one is that this limit is more explicitly given by

\begin{conjecture}\label{ConjAs2}
$$
 A(\Omega)\lim_{k\rightarrow +\infty} \frac{\mathfrak
 L_k(\Omega)}{k}=\lambda(\mathrm{Hexa}_1)\;.
$$
\end{conjecture}

\noindent The second conjecture says in particular that the limit
only depends on the area of $\Omega$, not on its geometry (provided
$\Omega$ is a regular domain).

\begin{figure}[ht]\label{F-1}
\begin{center}
\includegraphics[width=8cm]{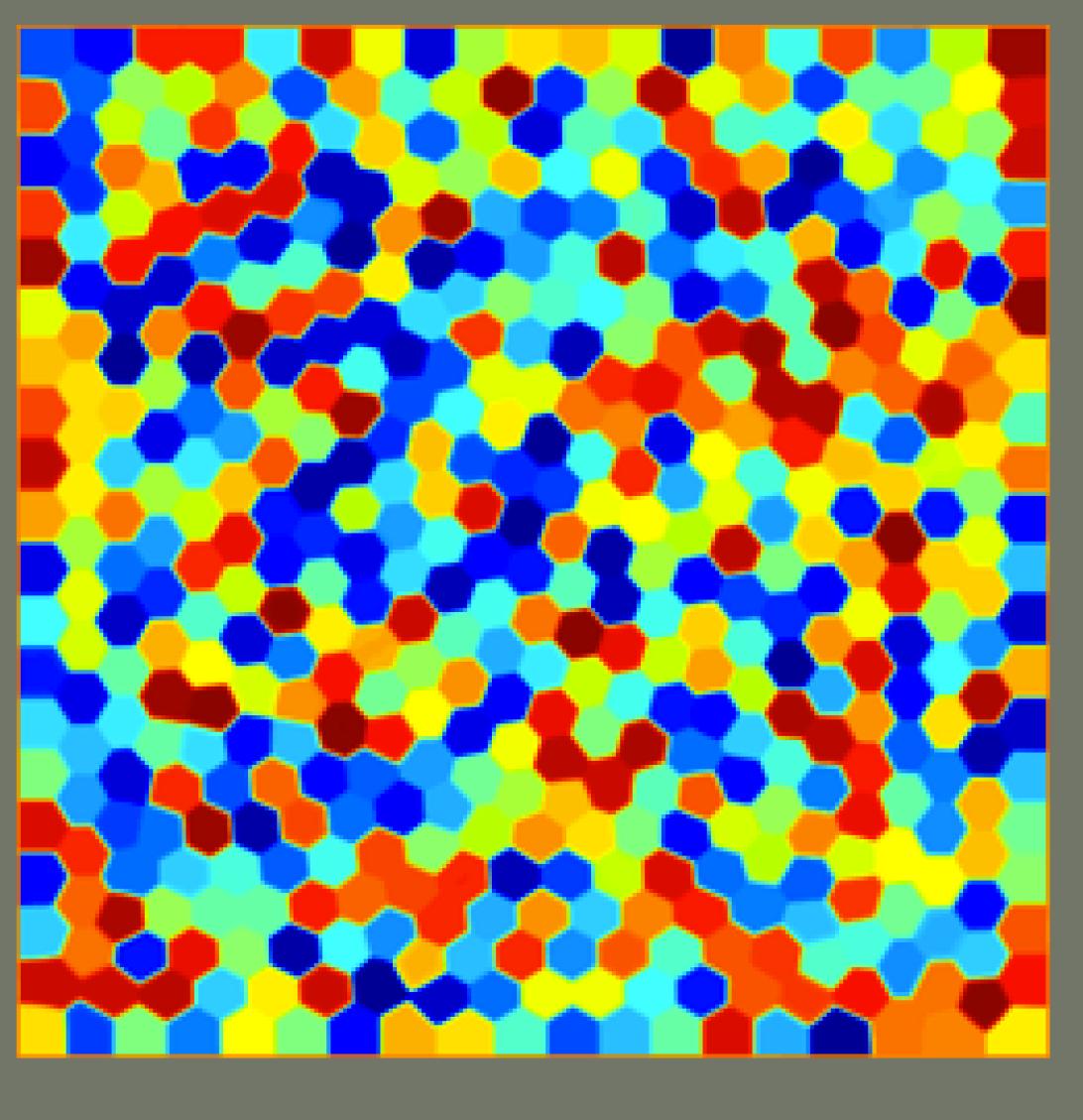}
\caption{Computations of  Bourdin-Bucur-Oudet for the periodic
square. (Minimization of the sum)}
\end{center}
\end{figure}

\medskip It is  explored numerically in \cite{BHV} why the second conjecture
looks reasonable. Note that Caffarelli and Lin \cite{CL1}, mention
Conjecture~\ref{ConjAs2} in relation with $\mathfrak
L_{k,1}(\Omega)$. From this point of view, the recent numerical
computations by Bourdin-Bucur-Oudet \cite{BBO} for the asymptotic
structure of the minimal partitions for $\mathfrak L_{k,1}(\Omega)$
are very enlightning, see Figure~\ref{F-1}. Remark~\ref{R-fg1} shows
that \eqref{E-glb} should be a strict inequality.\medskip

\subsection{Asymptotics of the length of the boundary set of\\  minimal
regular $k$-equipartitions for $k$ large.}

In this Section, we only consider bounded open domains $\Omega$ in $\rz^2$.
In this case, Conjecture~\ref{ConjAs2} leads to a natural ``hexagonal
conjecture'' for the length of the boundary set, namely

\begin{conjecture}\label{ConjAs3}
\begin{equation}\label{hcl}
\lim_{k\rightarrow +\infty}  ( P(\mathcal D_k)/\sqrt{k}) =
\frac{1}{2} \ell (\mathrm{Hexa}_1)  \sqrt{ A(\Omega)}\,,
\end{equation}
where $\ell (\mathrm{Hexa}_1)$ is the length of the boundary of the
hexagon of area $1$,
$$
\ell (\mathrm{Hexa}_1) = 2 \sqrt{2\sqrt{3}} = 2\, (12)^{\frac{1}{4}}\,.
$$
\end{conjecture}

For regular spectral equipartitions $\mathcal{D}$ of the domain $\Omega$,
inequalities \eqref{bg} and \eqref{E-fke} yield,
\begin{equation}\label{bga}
\liminf_{\sharp(\mathcal{D}) \to \infty} \frac{P(\mathcal{D})}{\sqrt{\sharp(\mathcal{D})}}
\ge \frac{\sqrt{\pi}}{2} \sqrt{A(\Omega)}.
\end{equation}
Assuming that $\chi(\Omega) \ge 0$, we have the uniform lower bound,
\begin{equation}\label{bga-1}
\frac{P(\mathcal{D})}{\sqrt{\sharp(\mathcal{D})}} \ge
\frac{\sqrt{\pi}}{2} \sqrt{A(\Omega)}.
\end{equation}

\noindent \textbf{Remark}. Assume that $\chi(\Omega) \ge 0$, and
that all the sub-domains $D_i$ in the regular equipartition
$\mathcal{D}$ satisfy $\chi(D_i) \ge 0$ as well. Then, Polya's
inequality \eqref{Polya} yields the sharper inequality,
\begin{equation}\label{bga-2}
\frac{P(\mathcal{D})}{\sqrt{\sharp(\mathcal{D})}} \ge
\frac{\mathbf{j}}{\sqrt{\pi}} \sqrt{A(\Omega)}.
\end{equation}\medskip

The following statement is a particular case of Theorem~1-B established by
T.C.~Hales \cite{Ha} in his proof of Lord Kelvin's honeycomb conjecture (see also \cite{Mor}).

\begin{theorem}
Let $\Omega$ be a relatively compact open set in $\mathbb{R}^2$, and
let $\mathcal{D}=\{D_i\}$  be a regular finite partition of $\Omega$. Then,
\begin{equation}
 P(\mathcal{D}) + \frac 12 \ell(\pa \Omega) \geq (12)^{\frac{1}{4}} \sum_{i=1}
 ^{\sharp(\mathcal{D})} \min \left(1,A(D_i)\right)\,.
\end{equation}
\end{theorem}

\noindent In order to optimize the use of this theorem, we consider $\min_i
A(D_i)$, and apply the theorem to a dilated partition. If we dilate
by $t$ the length is multiplied by $t$ and the area by $t^2$. So we
take  $ t=\left(\min_i A(D_i) \right)^{-\frac{1}{2}}$, and we obtain,

\begin{corollary}\label{C-hales}
For any regular partition $\mathcal{D}$ of a bounded open subset $\Omega$
of $\rz^2$,
\begin{equation}\label{honey1}
P(\mathcal D) + \frac{1}{2} \ell(\pa \Omega)  \geq (12)^{\frac{1}{4}}
 \, (\min_i A(D_i) )^{\frac{1}{2}} \, \sharp(\mathcal{D})\,.
\end{equation}
\end{corollary}

\begin{proposition}\label{P-ale}
Let $\Omega$ be a regular bounded domain in $\rz^2$.
\begin{enumerate}
    \item[(i)] For $k\ge 1$, let $\mathcal{D}_k$ be a minimal regular $k$-equipartition
    of $\Omega$. Then,
    \begin{equation}\label{asym}
    \liminf_{k\ar +\infty} \frac { P(\mathcal D_k)}{\sqrt{k}}  \geq
    (12)^{\frac{1}{4}} \left( \frac{\pi{\bf
    j}^2}{\lambda(\mathrm{Hexa}_1)} \right)^{\frac{1}{2}}\,
    A(\Omega)^{\frac{1}{2}}\,.
    \end{equation}
    \item[(ii)] If $\chi(\Omega) \ge 0$, then for any regular
    spectral equipartition, we have the universal estimate
    \begin{equation}\label{univ}
    P(\mathcal{D})+ \frac{1}{2} \ell(\pa \Omega)   \geq  12^{\frac{1}{8}}\,
    (\frac{\pi}{4})^{\frac{1}{4}} \,  A(\Omega)^{\frac{1}{2}}\, \left(\sharp(\mathcal{D})
    \right)^{\frac{1}{2}} \,.
    \end{equation}
\end{enumerate}
\end{proposition} \medskip

\noindent \textbf{Remarks}. \\
(a)~ Recall that $\lambda(\mathrm{Hexa}_1) \sim 18,5901$ and that
$\lambda (\mathrm{Disk}_1) = \pi {\bf j}^2 \sim 18,1680$. It follows
that $\left( \frac{\pi{\bf j}^2}{\lambda(\mathrm{Hexa}_1)}
\right)^{\frac{1}{2}}\sim 0,989\,$, so that the right-hand side of
\eqref{asym} is very close to the right-hand side of \eqref{hcl},
the hexagonal conjecture for the length.\\
 (b)~ Asymptotically, when
we consider minimal regular $k$-equipartition $\mathcal{D}_k$,
inequality \eqref{univ} is weaker than \eqref{asym} but it is
universal, and independent of the asymptotics of the energy of the
partition.\medskip

\noindent \textbf{Proof of the Proposition.}\\
  (i)~ Let $\mathcal{D}
=\{D_i\}$ be a regular equipartition of $\Omega$. Faber-Krahn's
inequality \eqref{E-fki} gives,
$$
\Lambda(\mathcal{D}) = \lambda (D_i) \geq \frac{\pi {\bf
j^2}}{A(D_i)}\,,
$$
hence, combining with \eqref{honey1},
\begin{equation*}\label{haa}
 P(\mathcal{D})+ \frac{1}{2} \ell(\pa \Omega)  \geq (12)^{\frac{1}{4}} \,
 (\pi {\bf j}^2)^{\frac{1}{2}} \frac{\sharp(\mathcal{D})}{\sqrt{\Lambda(\mathcal{D})}}\,.
\end{equation*}

\noindent Let $\mathcal{D}_k$ be a minimal regular $k$-equipartition of $\Omega$.
Applying \eqref{haa}, we obtain
\begin{equation*}
\frac{P(\mathcal{D}_k)+ \frac{1}{2} \ell(\pa \Omega)}{\sqrt{k}}   \geq (12)^{\frac{1}{4}}
\, (\pi {\bf j}^2)^{\frac{1}{2}} \, \big(\frac{\mathfrak L_k(\Omega)}{k}\big)^{-\frac{1}{2}}\,.
\end{equation*}

\noindent Using the upper bound for $\limsup_{k\ar +\infty}
\frac {\mathfrak L_k(\Omega)}{k}$ given by \eqref{E-asym-l}, we obtain
the following asymptotic inequality for the length of a minimal
regular $k$-equipartition,
\begin{equation*}
\liminf_{k\ar +\infty} \frac { P(\mathcal D_k)}{\sqrt{k}}  \geq
(12)^{\frac{1}{4}} \left( \frac{\pi{\bf
j}^2}{\lambda(\mathrm{Hexa}_1)} \right)^{\frac{1}{2}}\,
A(\Omega)^{\frac{1}{2}}\,.
\end{equation*}

(ii)~ Assume that $\chi(\Omega) \ge 0$. Indeed, by \eqref{bg}, we
have
$$
P(\mathcal{D}) \ge \frac{A(\Omega)}{2\mathbf{j}}
\sqrt{\Lambda(\mathcal{D})},
$$
and hence,
$$
\Lambda(\mathcal{D})^{-\frac{1}{2}} \ge \frac{A(\Omega)}{2
\mathbf{j} \, P(\mathcal{D})} \ge \frac{A(\Omega)}{2 \mathbf{j} \,
\big( P(\mathcal{D}) + \frac{1}{2}\ell(\partial\Omega)\big)}.
$$

\noindent The same proof as for Assertion (i) gives, for any regular
spectral equipartition,
\begin{equation*}
 P(\mathcal D)+ \frac{1}{2} \ell(\pa \Omega)  \geq (12)^{\frac{1}{4}} \,
 (\pi {\bf j}^2)^{\frac{1}{2}} \Lambda(\mathcal{D})^{-\frac{1}{2}}\,.
\end{equation*}
Inequality \eqref{univ} follows. \qed \medskip

\noindent \textbf{Remark}. Assume now that $\mathcal{D}(u_k)$ is the
nodal partition of some $k$-th eigenfunction $u_k$ of $H(\Omega)$.
Assume furthermore that $\chi(\Omega) \ge 0$. The same reasoning as
above gives,
\begin{equation*}
 P(\mathcal{D})+ \frac{1}{2} \ell(\pa \Omega)   \geq (12)^{\frac{1}{4}} \mu(u_k)
 \, (\pi {\bf j}^2)^{\frac 12} (\lambda_k(\Omega))^{-\frac{1}{2}}\,.
\end{equation*}
Hence
\begin{equation*}
\liminf_{k\ar +\infty} \frac { P\left(\mathcal{D}(u_k)\right)}{\sqrt{k} }
\geq (12)^{\frac 14} \left (\liminf_{k\ar +\infty}  \frac{\mu(u_k)}{k} \right)
\, ( {\bf j}^2/4)^{\frac 12} A(\Omega) ^{\frac 12}\,.
\end{equation*}
This inequality is less convincing, because we have no lower bound
for the right-hand side. Indeed, on the round sphere $\mathbb{S}^2$, any
eigenspace of the Laplacian contains eigenfunctions with only either $2$ or $3$
nodal domains, \cite{LEW}. \medskip

\noindent We may however go back to the initial inequality,
\begin{equation*}
P\left( \mathcal{D}(u_k)\right) \ge \frac{A(\Omega)}{2\mathbf{j}}
\sqrt{\lambda_k(\Omega)}.
\end{equation*}

\noindent Assuming Polya's conjecture \cite{Pol1} which says that for bounded
domains in $\rz^2$, $\lambda_k(\Omega) \ge \frac{4 \pi}{A(\Omega)}$,
for $k \ge 1$, we find that
\begin{equation*}
P\left( \mathcal{D}(u_k)\right) \ge \frac{\sqrt{\pi}}{\mathbf{j}} \, \sqrt{A(\Omega)}
\, \sqrt{k}\,,
\end{equation*}
which should be compared to \eqref{bga-1}.

\vspace{1cm}

%\bibliography{Bibliog/bib}
%\newpage

\end{document}